\documentstyle[12pt, amsfonts]{amsart}
\begin{document}
\def\R{{\mathbb R}}
\def\Z{{\mathbb Z}}
\def\C{{\mathbb C}}
\newcommand{\trace}{\rm trace}
\newcommand{\Ex}{{\mathbb{E}}}
\newcommand{\Prob}{{\mathbb{P}}}
\newcommand{\E}{{\cal E}}
\newcommand{\F}{{\cal F}}
\newtheorem{df}{Definition}
\newtheorem{theorem}{Theorem}
\newtheorem{lemma}{Lemma}
\newtheorem{pr}{Proposition}
\newtheorem{co}{Corollary}
\newtheorem{pb}{Problem}
\def\n{\nu}
\def\sign{\mbox{ sign }}
\def\a{\alpha}
\def\N{{\mathbb N}}
\def\A{{\cal A}}
\def\L{{\cal L}}
\def\X{{\cal X}}
\def\F{{\cal F}}
\def\c{\bar{c}}
\def\v{\nu}
\def\d{\delta}
\def\diam{\mbox{\rm dim}}
\def\vol{\mbox{\rm Vol}}
\def\b{\beta}
\def\t{\theta}
\def\l{\lambda}
\def\e{\varepsilon}
\def\colon{{:}\;}
\def\pf{\noindent {\bf Proof :  \  }}
\def\endpf{ \begin{flushright}
$ \Box $ \\
\end{flushright}}

\title[Inequalities for projections]
{Stability inequalities for projections of convex bodies}

\author{Alexander Koldobsky}

\address{Department of Mathematics\\ 
University of Missouri\\
Columbia, MO 65211}

\email{koldobskiya@@missouri.edu}

\begin{abstract}  The {projection function} 
$P_K$ of an origin-symmetric convex body $K$ in $\R^n$ is defined by
$P_K(\xi)=|K\vert {\xi^\bot}|,\ \xi\in S^{n-1},$
where $K\vert {\xi^\bot}$ is the projection of $K$ to the central hyperplaner $\xi^\bot$ perpendicular to $\xi$, and $|K|$ stands for volume of proper dimension.

We prove several stability and separation results for the projection function. For example,
if $D$ is a projection body in $\R^n$ which is in isotropic position up to a dilation,
and $K$ is any origin-symmetric convex body in $\R^n$ such that that there exists $\xi\in S^{n-1}$
with $P_K(\xi)>P_D(\xi),$ then
$$
\max_{\xi\in S^{n-1}} (P_K(\xi)-P_D(\xi)) \ge \frac c{\log^2n} (|K|^{\frac {n-1}n} -|D|^{\frac {n-1}n}),
$$
where $c$ is an absolute constant.

As a consequence, we prove a hyperplane inequality
$$
S(D) \le\  C \log^2n \max_{\xi\in S^{n-1}} S(D\vert\xi^\bot)\ |D|^{\frac 1n},
$$
where $D$ is a projection body in isotropic position,
up to a dilation, $S(D)$ is the surface area of $D,$ and $C$ is an absolute constant.
The proofs are based on the Fourier analytic approach to projections developed
in \cite{KRZ}.

\end{abstract}  
\maketitle

\section{Introduction}
Let $K$ be an origin-symmetric convex body in $\R^n.$ The {projection function} 
$P_K: S^{n-1}\to (0,\infty)$ of $K$ is defined for every
$\xi\in S^{n-1}$ as the $(n-1)$-dimensional volume of the orthogonal projection of $K$ to the central hyperplane 
$\xi^\bot$ perpendicular to $\xi.$ We write
$$P_K(\xi)=|K\vert {\xi^\bot}|,\qquad \forall \xi\in S^{n-1},$$
where $K\vert {\xi^\bot}$ is the projection of $K$ to $\xi^\bot$, and $|K|$ stands for volume of proper dimension.

The classical {\it uniqueness theorem} of Aleksandrov \cite{A} states that every origin-symmetric 
convex body is uniquely determined by the function $P_K;$ see also \cite[Theorem 3.3.6]{G}.
The corresponding  {\it volume comparison} question was posed by Shephard \cite{Sh} in 1964. 
Suppose that $K,D$ are origin-symmetric convex bodies in $\R^n,$ and 
$P_K(\xi)\le P_D(\xi)$ for every $\xi\in S^{n-1}.$ Does it necessarily follow that
$|K|\le |D| ?$ Shephard's problem was solved by Petty \cite{Pe} and Schneider \cite{S1},
independently, and the answer is affirmative only in dimension $n=2.$ Both solutions 
were based on a connection with projection bodies (see definition in Section \ref{stab-sep}), as follows.
If $D$ is a projection body, the answer to Shephard's question is affirmative for every $K.$ On the other hand,
if $K$ is not a projection body, one can construct $D$ giving together with $K$ a counterexample.
The final answer to the problem follows from the fact that only in dimension $n=2$ all origin-symmetric convex 
bodies are projection bodies.

For the affirmative cases in volume comparison problems, stability
and separation problems were proposed in \cite{K2}. Suppose that $\e>0,$ $K,L$ are origin-symmetric 
convex bodies in $\R^n,$ and $L$ is a projection body. {\it The stability problem} asks whether there exists a constant $c>0$ not dependent
on $K,L$ or $\e$ and such that the inequalities 
$P_K(\xi)\le P_L(\xi) + \e,$ for all $\xi\in S^{n-1},$
imply 
$|K|^{\frac{n-1}n}\le |L|^{\frac{n-1}n} + c\e.$
{\it The separation problem} asks whether inequalities 
$P_K(\xi)\le P_L(\xi) - \e,$ for all $\xi\in S^{n-1},$
imply 
$|K|^{\frac{n-1}n}\le |L|^{\frac{n-1}n} - c\e,$ again
with a constant $c>0$ not dependent on $K,L,\e.$

Separation was proved in \cite{K2} (see 
also \cite{K3}, where the constant is written precisely)
with the best possible constant 
$$c=c_n=\frac{|B_2^n|^{\frac{n-1}n}}{|B_2^{n-1}|},\ c_n\in (\frac 1{\sqrt{e}}, 1),$$ 
where $B_2^n$ 
is the unit Euclidean ball in $\R^n.$ Stability 
was proved in \cite{K2} with constants dependent on the bodies, using rough estimates for $M$ and $M^*$
parameters of convex bodies. In this article we prove stability with constants not dependent on the bodies,
but under an additional assumption that $L$ is in isotropic position, up to a dilation; see Proposition \ref{stabproj}. 
To do this, we use a recent 
$M^*$-estimate of Milman \cite{M}. We also prove that, if the projection body condition is dropped,
one can get results going in the direction opposite to stability and separation.
Namely, we give examples of origin-symmetric convex  bodies $K,L$ such that 
$P_K\le P_L+\e,$ but $|K|^{\frac{n-1}n}\ge |L|^{\frac{n-1}n} + c\e,$ and also examples of $K,L$
where $P_K\le P_L-\e,$ but $|K|^{\frac{n-1}n}\ge |L|^{\frac{n-1}n} - c\e,$ with $c$ not dependent on the bodies 
or small enough $\e.$ In some sense, these results provide a quantitative version of the solution to
Shephard's problem.

Stability and separation immediately imply what we call {\it volume difference inequalities}. In fact, if stability holds and there exists $\xi\in S^{n-1}$ such that $P_K(\xi)>P_L(\xi),$ then we put 
$\e=\max_{\xi\in S^{n-1}} \left(P_K(\xi)-P_L(\xi)\right),$
and get
$$|K|^{\frac{n-1}n} - |L|^{\frac{n-1}n} \le c \max_{\xi\in S^{n-1}} \left(P_K(\xi)-P_L(\xi)\right).$$
Similarly, if separation holds, we get
$$|L|^{\frac{n-1}n} - |K|^{\frac{n-1}n} \ge c \min_{\xi\in S^{n-1}} \left(P_L(\xi)-P_K(\xi)\right).$$
We provide such an inequality for the stability result mentioned above; see Theorem \ref{th-max-proj}.

Volume difference inequalities lead to {\it hyperplane inequalities} for surface area of projection
bodies. It was proved in \cite{K3} that
if $D$ is a projection body in $\R^n,$ then 
\begin{equation}\label{surface}
S(D) \ge c \min_{\xi\in S^{n-1}} S(D\vert\xi^\bot)\ |D|^{\frac 1n},
\end{equation}
where $S(D)$ is the surface area of $D,$ and $c$ is an absolute constant.
In this paper we prove an inequality that complements (\ref{surface}).
If $D$ is a projection body in isotropic position, up to a dilation, then
$$
S(D) \le\  C \log^2n \max_{\xi\in S^{n-1}} S(D\vert\xi^\bot)\ |D|^{\frac 1n},
$$
where $C$ is an absolute constant; see Theorem \ref{th-hyper-max} below.

Finally, we show that if the condition that $D$ is a projection body is removed,
volume difference inequalities can go in the opposite direction; see Theorems \ref{counter-max}
and \ref{counter-min}.

We extensively use the Fourier analytic
approach to projections of convex bodies developed in \cite{KRZ}; see also \cite[Chapter 8]{K1}.

\section{Stability theorems} \label{stab-sep}

We need several definitions from convex geometry. We refer the reader
to \cite{S2} for details.

The {\it support function} of a convex body $K$ in $\R^n$ is defined by
$$h_K(x) = \max_{\{\xi\in \R^n:\|\xi\|_K=1\}} (x,\xi),\quad x\in \R^n.$$ 
If $K$ is origin-symmetric, then $h_K$ is a norm on $\R^n.$

The {\it surface area measure} $S(K, \cdot)$ of a convex body $K$ in 
$\R^n$ is a measure on $S^{n-1}$ defined as follows. For every Borel set $E \subset S^{n-1},$ 
$S(K,E)$ is equal to Lebesgue measure of the part of the boundary of $K$
where normal vectors belong to $E.$ 
We usually consider bodies with absolutely continuous surface area measures.
A convex body $K$ is said to have the {\it curvature function} 
$ f_K: S^{n-1} \to \R$
if its surface area measure $S(K, \cdot)$ is absolutely 
continuous with respect to Lebesgue measure $\sigma_{n-1}$ on 
$S^{n-1}$, and
$$
\frac{d S(K, \cdot)}{d \sigma_{n-1}}=f_K \in L_1(S^{n-1}),
$$
so $f_K$ is the density of $S(K,\cdot).$

The volume of a body can be expressed in terms of its support function and 
curvature function:
\begin{equation}\label{vol-proj}
|K| = \frac 1n \int_{S^{n-1}}h_K(x) dS(K,x) = \frac 1n \int_{S^{n-1}}h_K(x) f_K(x)\ dx,
\end{equation}
with the latter equality if $f_K$ exists.

If $K$ and $L$ are two convex bodies in $\R^n,$ the {mixed volume} $V_1(K,L)$
is equal to  
$$V_1(K,L)= \frac{1}{n} \lim_{\e\to +0}
\frac{|K+\epsilon L|-|K|}{\e}.$$
The first Minkowski inequality (see for example \cite[p.23]{K1}) asserts that
for any convex bodies $K,L$ in $\R^n,$ 
\begin{equation} \label{firstmink}
V_1(K,L) \ge |K|^{\frac{n-1}n} |L|^{\frac1n}.
\end{equation}
Mixed volume can also be expressed in terms of the support and
curvature functions:

\begin{equation}\label{mixvol-proj}
V_1(K,L) = \frac 1n \int_{S^{n-1}}h_L(x) dS(K,x)= \frac 1n \int_{S^{n-1}}h_L(x) f_K(x)\ dx.
\end{equation}

Let $K$ be an origin-symmetric convex body in $\R^n.$ The {\it
projection body} $\Pi K$ of $K$ is defined by
\begin{equation} \label{def:proj}
h_{\Pi K}(\theta) = |K\vert {\theta^{\perp}}|,\qquad \forall \theta\in S^{n-1}.
\end{equation}

If $L$ is the projection body of some convex body, we simply say 
that $L$ is a projection body. We refer the reader to \cite[Chapter 8]{K1} for necessary information
about projection bodies. We just mention here that, by a result of Bolker \cite{Bl}, an origin-symmetric
convex body is a projection body if and only if its polar body is the unit ball of a subspace of $L_1.$
The unit balls of the spaces $\ell_p^n,\ p\ge 2$ are projection bodies, while the unit balls of
$\ell_p^n,\ p<2,\ n\ge 3$ are not. There are other examples of bodies in $\R^n,\ n\ge 3$ that are not
projection bodies; see \cite{S2} or \cite{K1}.

The classes of projection bodies $K$ for which the functions $h_K$ and $f_K$ are infinitely smooth
are dense in the class of all projection bodies in the Hausdorff metric; see \cite[p.151]{S2} and \cite{GZ}.
The class of projection bodies for which the curvature function is strictly positive is also
dense in the class of all projection bodies; see \cite{GZ} or \cite[p.158, p. 161]{K1}.

We say that a body $K$ is in {isotropic position} if $|K|=1$ and there exists a constant $L_K>0$
such that
$$\int_K (x,\xi)^2 dx = L_K^2,\qquad \forall \xi\in S^{n-1}.$$
For every origin-symmetric convex body $K$ in $\R^n$ there exists $T\in GL_n$ such that
$TK$ is in isotropic position. The constant $L_K$ is called the isotropic
constant of $K.$ It is known that $L_K\ge L_{B_2^n}$ for every symmetric convex body $K$
in $\R^n;$ see \cite[p.93]{MP}. The question of whether $L_K$ is bounded from above by an absolute
constant is the matter of the well-known and still open slicing problem. The best known
estimate $L_K\le O(n^{1/4})$ is due to Klartag \cite{Kl}, who improved an earlier estimate of
Bourgain \cite{Bo}. We refer the reader to \cite{BGVV} for these results and more about the
isotropic position and the slicing problem.

We use the Fourier approach to projections of convex bodies developed in \cite{KRZ}; see also \cite[Chapter 8]{K1}.
We consider Schwartz distributions, i.e. continuous functionals on the space ${\cal{S}}(\R^n)$
of rapidly decreasing infinitely differentiable functions on $\R^n$. 
The Fourier transform of a distribution $f$ is defined by $\langle\hat{f}, \phi\rangle= \langle f, \hat{\phi} \rangle$ for
every test function $\phi \in {\cal{S}}(\R^n).$ For any even distribution $f$, we have $(\hat{f})^\wedge
= (2\pi)^n f$.

For $f\in C^\infty(S^{n-1})$ and $p=1$ or $p=-n-1,$ we denote by 
$$(f\cdot r^{p})(x) = f(x/|x|_2) |x|_2^{p}$$
the extension of $f$ to a homogeneous function of degree $p$ on $\R^n.$
By  \cite[Lemma 3.16]{K1}, there exists $g\in C^\infty(S^{n-1})$ such that
\begin{equation} \label{homogeneous}
(f\cdot r^{p})^\wedge = g\cdot r^{-n-p}.
\end{equation}

In particular, if the support function $h_K$ is infinitely smooth on $S^{n-1},$ then
the Fourier transform of $h_K\cdot r$ is an infinitely smooth function $g$ on the sphere
extended to a homogeneous function of degree $-n-1$ on $\R^n\setminus \{0\}.$
In this case we simply write $(h_K\cdot r)^\wedge(\theta)$ for $\theta\in S^{n-1},$
meaning the function $g.$ It was proved in \cite{KRZ} (see also \cite[Theorem 8.6]{K1})
that an origin-symmetric convex body $K,$ for which $h_K$ is infinitely smooth,
is a projection body if and only if
\begin{equation} \label{projbody}
(h_K\cdot r)^\wedge(\theta)\le 0,\qquad \forall \theta\in S^{n-1}.
\end{equation}

It was also proved in \cite{KRZ} (see also \cite[Theorem 8.2]{K1}) that if an origin-symmetric  body $K$
has curvature function $f_K$, then
\begin{equation}\label{proj-curv}
(f_K\cdot r^{-n-1})^\wedge(\theta)= -\pi |K\vert_{\theta^\bot}|=-\pi P_K(\theta),\qquad \forall \theta\in S^{n-1}.
\end{equation}

The following version of Parseval's formula was proved in \cite{KRZ} (see also \cite[Lemma 8.8]{K1}).
If $K,L$ are origin-symmetric convex bodies, $K$ has infinitely smooth support function and $L$ has
infinitely smooth curvature function, then 
\begin{equation} \label{pars-proj}
\int_{S^{n-1}} (h_K\cdot r)^\wedge (\xi) (f_L\cdot r^{-n-1})^\wedge(\xi)\ d\xi =
(2\pi)^n \int_{S^{n-1}} h_K(x) f_L(x)\ dx.
\end{equation}

\begin{lemma}\label{remain3} Let $K$ be an origin-symmetric convex body in $\R^n$ such that
the support function $h_K$ is infinitely smooth.
Then
$$\int_{S^{n-1}} (h_K\cdot r)^\wedge(\xi) d\xi \le -\frac{(2\pi)^n n}{\pi} c_n |K|^{1/n}.$$
Recall that $c_n=|B_2^n|^{\frac{n-1}n}/|B_2^{n-1}|\in (\sqrt{\frac 1{e}},1).$
\end{lemma}
\pf The curvature function of the unit Euclidean ball $B_2^n$ is constant, $f_2\equiv 1.$
By (\ref{proj-curv}), 
$$(f_2\cdot r^{-n-1})^\wedge(\xi)= -\pi |B_2^{n-1}|,\qquad \forall \xi\in S^{n-1}.$$
By (\ref{mixvol-proj}), (\ref{firstmink}) and (\ref{pars-proj}),
$$\int_{S^{n-1}} (h_K\cdot r)^\wedge(\xi) d\xi= -\frac 1{\pi |B_2^{n-1}|} 
\int_{S^{n-1}} (h_K\cdot r)^\wedge(\xi) (f_2\cdot r^{-n-1})^\wedge(\xi) d\xi$$
$$= -\frac {(2\pi)^n}{\pi |B_2^{n-1}|} \int_{S^{n-1}} h_K(x) f_2(x) dx \le
-\frac {(2\pi)^n n |K|^{1/n} |B_2^n|^{(n-1)/n}}{\pi |B_2^{n-1}|}. \qed$$

\begin{lemma} \label{remain4}  Let $K$ be an origin-symmetric convex body in $\R^n$
such that the support function $h_K$ is infinitely smooth and $K$ is a dilate of an isotropic body. 
Then
$$\int_{S^{n-1}} (h_K\cdot r)^\wedge(\xi) d\xi \ge -C (2\pi)^n n \log^2(1+n) L_K |K|^{1/n},$$
where $C$ is an absolute constant.
\end{lemma}
\pf By the same argument as in Lemma \ref{remain3} (recall that $f_2\equiv 1)$ we have
$$\int_{S^{n-1}} (h_K\cdot r)^\wedge(\xi) d\xi=  -\frac {(2\pi)^n}{\pi |B_2^{n-1}|} \int_{S^{n-1}} h_K(x)dx.$$
Now use the following estimate of E.Milman \cite[Theorem 1.1]{M} 
$$\frac 1{|S^{n-1}|} \int_{S^{n-1}} h_K(x) dx \le C_1\sqrt{n} \log^2(1+n) L_K |K|^{1/n},$$
where $C_1$ is an absolute constant,
and note that $|S^{n-1}|/|B_2^{n-1}|\sim \sqrt{n}$ to get the result. \qed
\bigbreak
We now prove stability in the affirmative direction of Shephard's problem under the additional
condition that the body $D$ is a dilate of an isotropic body.
\begin{pr} \label{stabproj} Suppose that $\e>0$,  $K$ and $D$ are origin-symmetric
convex bodies in $\R^n,$ and $D$ is a projection body which is a dilate of an isotropic body.  
If for every $\xi\in S^{n-1}$
\begin{eqnarray}\label{proj11}
P_K(\xi)\le P_D(\xi) +\e,
\end{eqnarray}
then
$$|K|^{\frac{n-1}n}  \le |D|^{\frac{n-1}n} + C \e \log^2(1+n) L_D,$$
where $C$ is an absolute constant.
\end{pr}

\pf  By approximation (\cite[Th. 3.3.1]{S2} and \cite[Section 5]{GZ}), we can assume 
that $D$ has infinitely smooth support function and both $D$ and $K$ have
infinitely smooth curvature functions. By (\ref{proj-curv}), the condition (\ref{proj11})
can be written as
\begin{equation} \label{fourier-proj}
-\frac 1{\pi} (f_K\cdot r^{-n-1})^\wedge(\xi) \le  -\frac 1{\pi} (f_D\cdot r^{-n-1})^\wedge(\xi) + \e, \qquad \forall \xi\in S^{n-1}.
\end{equation}
By (\ref{projbody}),
$(h_D\cdot r)^\wedge \le 0$ on the sphere $S^{n-1}.$ Therefore, integrating (\ref{fourier-proj})
with respect to a negative density, we get
$$\int_{S^{n-1}} (h_D\cdot r)^\wedge(\xi) (f_D\cdot r^{-n-1})^\wedge(\xi)\ d\xi$$$$ \ge 
\int_{S^{n-1}} (h_D\cdot r)^\wedge(\xi) (f_K\cdot r^{-n-1})^\wedge(\xi)\ d\xi 
+ \pi\e \int_{S^{n-1}} (h_D\cdot r)^\wedge(\xi)\ d\xi.$$
Using this, (\ref{vol-proj}), (\ref{mixvol-proj}), Parseval's formula (\ref{pars-proj}) 
and the first Minkowski inequality (\ref{firstmink}), 
$$ (2\pi)^n n |D| = (2\pi)^n \int_{S^{n-1}} h_D(x) f_D(x)\  dx $$
$$
\ge (2\pi)^n \int_{S^{n-1}} h_D(x) f_K(x)\  dx + \pi\e\int_{S^{n-1}}(h_D\cdot r)^\wedge(\xi)\ d\xi.
$$
$$\ge  (2\pi)^n n |D|^{\frac 1n} |K|^{\frac {n-1}n}+ \pi\e\int_{S^{n-1}}(h_D\cdot r)^\wedge(\xi)\ d\xi.$$
The result follows from Lemma \ref{remain4}. \qed

Now we show that if the projection body condition is dropped, the result
may go in the opposite direction.

\begin{pr}\label{stabprojcounter} Suppose that $K$ is an origin-symmetric convex body in $\R^n$, 
which is not a projection body. Suppose also that $h_K, f_K\in C^\infty(S^{n-1}),$
and $f_K$ is strictly positive on $S^{n-1}.$ 
Then for small enough
$\e>0$ there exists an origin-symmetric convex body $D$ in $\R^n$ so that
$$P_D(\theta)\le P_K(\theta) \le P_D(\xi) + \e,\qquad \forall \theta\in S^{n-1},$$
but
$$|K|^{\frac{n-1}n}  \ge |D|^{\frac{n-1}n} + c_n \e.$$
\end{pr}

\pf By (\ref{projbody}), since $K$ is not a projection body,  there exists a symmetric open set $\Omega\subset S^{n-1},$
where $(h_K\cdot r)^\wedge>0.$  Let $v$ be an even infinitely smooth 
non-negative function supported in
$\Omega.$ Extend $v$ to a homogeneous function $v\cdot r$ of degree 1 on $\R^n.$ The
Fourier transform of $v\cdot r$ is a homogeneous of degree $-n-1$ function $g\cdot r^{-n-1},$
where $g$ is an infinitely smooth function on the sphere; recall (\ref{homogeneous}).

Define an even function $h$ on the sphere $S^{n-1}$ by
$$f_K=h+\delta g + \frac{\e}{|B_2^{n-1}|}.$$
Choose $\e,\delta>0$ small enough so that $h>0$ everywhere on $S^{n-1}$  
(recall that $f_K>0$). By the Minkowski existence theorem (see, for example, \cite[Section 7.1]{S2}), 
there exists an origin-symmetric convex 
body $D$ in $\R^n,$ whose curvature function $f_D=h.$ Extend the functions in the definition of $h$
to even homogeneous of degree $-n-1$ functions on $\R^n:$
\begin{equation}\label{def-h}
f_K\cdot r^{-n-1} = f_D\cdot r^{-n-1} + \delta g\cdot r^{-n-1} + \frac{\e}{|B_2^{n-1}|}\cdot r^{-n-1}.
\end{equation}
By (\ref{proj-curv}), 
$(r^{-n-1})^\wedge(\theta) = (f_{2}\cdot r^{-n-1})^\wedge(\theta) = -\pi |B_2^{n-1}|$
for every $\theta\in S^{n-1}$ (here $f_2\equiv 1$ is the curvature function of the unit Euclidean ball).

Taking the Fourier transform of both sides of (\ref{def-h}) and again using (\ref{proj-curv}), we get
\begin{equation}\label{def-D}
(f_K\cdot r^{-n-1})^\wedge(\theta) = (f_D\cdot r^{-n-1})^\wedge(\theta) + 
(2\pi)^n \delta\ v(\theta)- \pi\e,\qquad \forall \theta\in S^{n-1},
\end{equation}
and
$$-\pi P_K(\theta)= -\pi P_D(\theta) + (2\pi)^n \delta\ v(\theta) - \pi\e,\qquad \forall \theta\in S^{n-1}.$$
Since $v\ge 0,$ the latter implies
$$P_K(\theta) \le P_D(\theta) + \e,\qquad \forall \theta\in S^{n-1}.$$
Also, choosing $\delta$ small enough, we can assure that $P_K(\theta)\ge P_D(\theta)$ for every $\theta\in S^{n-1}.$

On the other hand, multiplying (\ref{def-D}) by $(h_K\cdot r)^\wedge(\theta)$, integrating over the sphere, and using
the fact that $v\ge 0$ is supported in $\Omega,$ where $(h_K\cdot r)^\wedge > 0,$  we get
$$ \int_{S^{n-1}} (h_K\cdot r)^\wedge(\theta)(f_K\cdot r^{-n-1})^\wedge(\theta) d\theta$$
$$= \int_{S^{n-1}} (h_K\cdot r)^\wedge(\theta)(f_D\cdot r^{-n-1})^\wedge(\theta) d\theta$$$$ + 
(2\pi)^n \delta \int_{S^{n-1}} v(\theta) (h_K\cdot r)^\wedge(\theta) d\theta -
\pi\e \int_{S^{n-1}} (h_K\cdot r)^\wedge(\theta)d\theta$$
$$\ge \int_{S^{n-1}} (h_K\cdot r)^\wedge(\theta)(f_D\cdot r^{-n-1})^\wedge(\theta)d\theta-
\pi\e \int_{S^{n-1}} \widehat{h_K}(\theta)d\theta.$$
Now, by
Parseval's formula (\ref{pars-proj}), (\ref{vol-proj}), (\ref{mixvol-proj}) and the first Minkowski inequality,
$$(2\pi)^n n |K| \ge (2\pi)^n \int_{S^{n-1}} h_K(\theta)f_D(\theta)d\theta -
\pi\e \int_{S^{n-1}}  (h_K\cdot r)^\wedge(\theta)d\theta$$
$$\ge (2\pi)^n n |K|^{\frac 1n}|D|^{\frac{n-1}n} -\pi\e \int_{S^{n-1}}  (h_K\cdot r)^\wedge(\theta)d\theta.$$
The result follows from Lemma \ref{remain3}. \qed
\smallbreak

The following separation result will be used to prove Theorem \ref{counter-min}.
\begin{pr}\label{sepprojcounter} Suppose that $K$ is an origin-symmetric convex body in $\R^n$
with strictly positive curvature which is not a projection body and is a dilate of an isotropic body.
Also suppose that the support and curvature functions of $K$ are infinitely smooth. Then for small enough
$\e>0$ there exists an origin-symmetric convex body $D$ in $\R^n$ so that 
$$P_K(\theta) \le P_D(\theta) - \e,\qquad \forall \theta\in S^{n-1},$$
but
$$|D|^{\frac{n-1}n}  \le |K|^{\frac{n-1}n} + C\log^2(1+n)L_K \e,$$
where $C$ is an absolute constant.
\end{pr}

\pf The proof follows the steps of the proof of Theorem \ref{stabprojcounter}.
Define the functions $v$ and $g$ in the same way.
Then define a body $D$ by
$$f_D=f_K - \delta g + \frac{\e}{|B_2^{n-1}|}.$$
At the very end use Lemma \ref{remain4} instead of Lemma \ref{remain3}. \qed

\section{Volume difference and hyperplane inequalities}

In this section we apply stability theorems to prove our main results. We start with the volume difference
inequality of Theorem \ref{th-max-proj}.
\medbreak
\begin{theorem} \label{th-max-proj} Let $D$ be a projection body in $\R^n$ in isotropic position up to a dilation,
and let $K$ be any origin-symmetric convex body in $\R^n.$ Suppose that there exists $\xi\in S^{n-1}$
so that $P_K(\xi)>P_D(\xi).$ Then
\begin{equation}\label{max-proj}
\max_{\xi\in S^{n-1}} (P_K(\xi)-P_D(\xi)) \ge \frac c{\log^2n} (|K|^{\frac {n-1}n} -|D|^{\frac {n-1}n}),
\end{equation}
where $c$ is an absolute constant.
\end{theorem}

\pf
Let $\e=\max_{\xi\in S^{n-1}} (P_K(\xi)-P_D(\xi)).$
By the condition of Theorem \ref{th-max-proj}, $\e>0.$ Now we can apply Proposition \ref{stabproj}
to $K,D, \e.$ We get 
$$|K|^{\frac{n-1}n}  \le |D|^{\frac{n-1}n} + C\log^2(1+n) L_D \max_{\xi\in S^{n-1}} (P_K(\xi)-P_D(\xi)),$$
where $C$ is an absolute constant. The result follows from the fact that isotropic constants of projection
bodies (zonoids) are uniformly bounded by an absolute constant; see \cite[p.96]{MP}. \qed
\smallbreak

``Differentiating" the inequality of Theorem \ref{th-max-proj}, we prove a hyperplane
inequality for the surface area of projection bodies.
\medbreak
\begin{theorem} \label{th-hyper-max}
Suppose that $D$ is a projection body in $\R^n$ in isotropic position up to a dilation (see definition
in Section \ref{stab-sep}). Then
\begin{equation}\label{hyper-max}
S(D) \le\  C \log^2n \max_{\xi\in S^{n-1}} S(D\vert\xi^\bot)\ |D|^{\frac 1n},
\end{equation}
where $C$ is an absolute constant. 
\end{theorem}

\pf
The surface area of $D$ can be computed as
$$S(D) = \lim_{\e \to +0} \frac {\left|D+\e B_2^n\right| - \left|D\right|}{\e}.$$
The inequality (\ref{max-proj}) with the bodies
$K=D+\e B_2^n$ and $D$ implies
\begin{equation} \label{surflimit}
\frac {|D+\e B_2^n|^{\frac {n-1}n} - |D|^{\frac {n-1}n}}{\e} \le
C\log^2n \max_{\xi\in S^{n-1}} \frac {|(D\vert \xi^\bot)+\e B_2^{n-1}| - 
|D\vert \xi^\bot|}{\e},
\end{equation}
where $C$ is an absolute constant.

By the Minkowski theorem on mixed volumes (\cite[Theorem 5.1.6]{S2} or \cite[Theorem A.3.1]{G}),
\begin{equation} \label{quer}
\frac{|(D\vert \xi^\bot)+ \e B_2^{n-1}|-|D\vert \xi^\bot|}{\e} = 
\sum_{i=1}^{n-1} {n-1 \choose i} W_i(D\vert \xi^\bot) \e^{i-1},
\end{equation}
where $W_i$ are quermassintegrals. The function $\xi\mapsto D\vert \xi^\bot$ is continuous 
from $S^{n-1}$ to the class of origin-symmetric convex sets equiped with the Hausdorff metric,
and $W_i$'s are also continuous with respect to this metric (see \cite[p.275]{S2}), so the 
functions $\xi\mapsto W_i(D\vert \xi^\bot)$ are continuous and, hence, bounded on the sphere.
This implies that the left-hand side of (\ref{quer}) converges to $S(D\vert \xi^\bot),$ as $\e\to 0,$
uniformly with respect to $\xi.$ The latter allows to switch the limit and maximum  
in the right-hand side of (\ref{surflimit}), as $\e\to 0$.
Sending $\e$ to zero in (\ref{surflimit}), we get
$$\frac {n-1}n |D|^{-1/n} S(D) \le C\log^2n  \max_{\xi \in S^{n-1}} S(D\vert \xi^\bot).$$
\endpf
\medbreak
Theorems \ref{counter-max} and \ref{counter-min} follow from Propositions \ref{stabprojcounter} 
and \ref{sepprojcounter} by putting 
$\e=\max_{\xi\in S^{n-1}} (P_K(\xi)-P_D(\xi))$
and 
$\e=\min_{\xi\in S^{n-1}} (P_D(\xi)-P_K(\xi)),$
correspondingly.

\begin{theorem}\label{counter-max}
Suppose that $K$ is an origin-symmetric convex body in $\R^n$ with strictly positive curvature
that is not a projection body. Then there exists an origin-symmetric convex body $D$
in $\R^n$ so that $P_K(\xi)\ge P_D(\xi)$ for all $\xi\in S^{n-1}$ and
$$\max_{\xi\in S^{n-1}} (P_K(\xi)-P_D(\xi))\le \frac 1{c_n} (|K|^{\frac {n-1}n}-|D|^{\frac {n-1}n}).$$
\end{theorem}

\begin{theorem} \label{counter-min} Suppose that $K$ is an origin-symmetric convex body in $\R^n$ 
that is not a projection body and is in isotropic position up to a dilation, with isotropic constant $L_K.$
Then there exists an origin-symmetric convex body $D$ in $\R^n$ so that $P_D(\xi)\ge P_K(\xi)$ for
all $\xi\in S^{n-1}$ and
$$\min_{\xi\in S^{n-1}} (P_D(\xi)-P_K(\xi)) \ge \frac c{L_K\log^2n} (|D|^{\frac {n-1}n}-|K|^{\frac {n-1}n}),$$
where $c$ is an absolute constant.
\end{theorem}

\noindent {\bf Remark.} Putting $D=\beta B_2^n$ in (\ref{max-proj}) and sending $\beta \to 0,$ we get a hyperplane inequality for volume
$$\max_{\xi\in S^{n-1}} P_K(\xi) \ge \frac c{\log^2n} |K|^{\frac {n-1}n}.$$
However, a stronger inequality
$$|L|^{\frac {n-1}n} \le c_n \max_{\xi\in S^{n-1}} |L\vert \xi^\bot|$$
holds for all origin-symmetric convex bodies and follows from the Cauchy projection formula for the surface area 
and the classical isoperimetric inequality; see \cite[p. 363]{G}.

It can be deduced directly from the solution to Shephard's problem (see \cite[Corollary 9.3.4]{G}) that,
if $L$ is a projection body in $\R^n$, then
\begin{equation} \label{hyper-proj}
|L|^{\frac {n-1}n} \ge c_n \min_{\xi\in S^{n-1}} |L\vert \xi^\bot|.
\end{equation}
Recall that $c_n > 1/\sqrt{e}.$ For general symmetric convex bodies, Ball \cite{Ba} proved that 
$c_n$ may and has to be replaced in (\ref{hyper-proj}) by
 $c/\sqrt{n},$ where $c$ is an absolute constant. 

\bigbreak
{\bf Acknowledgement.} I wish to thank the US National Science Foundation for support through 
grant DMS-1265155. Part of the work was done during my stay at the Max Planck Institute for
Mathematics in Bonn in Spring 2015.


\begin{thebibliography}{99}

\bibitem[A]{A} A.~D.~Aleksandrov, {\it On the surface area function of a convex body}, Mat. Sb. 
{\bf 6} (1939), 167--174.

\bibitem[Ba]{Ba} {K.~Ball}, \textit{Shadows of convex bodies}, Trans. Amer. Math. Soc.
\textbf{327}  (1991), 891--901. 

\bibitem[Bl]{Bl} E.~D.~Bolker, {\it A class of convex bodies,} Trans. Amer. Math. Soc. {\bf 145} (1969),
323--345.

\bibitem[Bo]{Bo} {J.~ Bourgain}, \textit{On the distribution of polynomials on high-dimensional 
convex sets},  Geometric aspects of functional analysis, Israel seminar (1989--90), Lecture Notes in Math.
\textbf{1469} Springer, Berlin, 1991, 127--137. 

\bibitem[BGVV]{BGVV} {S.~Brazitikos, A.~Giannopoulos, P.~Valettas and B.~Vritsiou},
{\it Geometry of isotropic log-concave measures}, Amer. Math. Soc., Providence, RI, 2014.

\bibitem[G]{G} { R.~J.~Gardner}, \textit{Geometric tomography}, Second edition,
Cambridge University Press, Cambridge, 2006. 

\bibitem[GZ]{GZ} {E.~Grinberg and Gaoyong Zhang,}{\it Convolutions, transforms, and convex bodies,} 
Proc. London. Math. Soc. (3) {\bf 78} (1999), 77--115.

\bibitem[Kl]{Kl} B.~Klartag, {\it On convex perturbations with a bounded isotropic constant}, 
Geom. Funct. Anal.  \textbf{16} (2006), 1274--1290. 

\bibitem[K1]{K1} {A.~Koldobsky}, \textit{Fourier analysis in convex geometry},
Amer. Math. Soc., Providence RI, 2005.

\bibitem[K2]{K2} A.~Koldobsky, {\it Stability in the Busemann-Petty and Shephard problems,} 
Adv. Math. {\bf 228} (2011), 2145--2161.

\bibitem[K3]{K3} A.~Koldobsky, {\it Stability and separation in volume comparison problems,} 
Math. Model. Nat. Phenom.  {\bf 8} (2013), 156--169.

\bibitem[KRZ]{KRZ} A.~Koldobsky, D.~Ryabogin and A.~Zvavitch, {\it Projections of convex
bodies and the Fourier transform}, Israel J. Math. {\bf 139} (2004), 361--380.

\bibitem[M]{M} E.~Milman, {\it On the mean-width of isotropic convex bodies and their
associated $L_p$-centroid bodies}, arxiv:1402.0209v1. 

\bibitem[MP]{MP} {V.~Milman and A.~Pajor}, {\em Isotropic position
and inertia ellipsoids and zonoids of the unit ball of 
a normed $n$-dimensional space}, in: Geometric Aspects
of Functional Analysis, ed. by J.~Lindenstrauss and V.~Milman,
Lecture Notes in Mathematics {\bf 1376}, Springer, Heidelberg, 1989, pp.~64--104.

\bibitem[Pe]{Pe} C.~M.~Petty, {\it Projection bodies}, Proc. Coll. Convexity (Copenhagen 1965),
Kobenhavns Univ. Mat. Inst., 234--241.

\bibitem[S1]{S1} R.~Schneider, {\it Zu einem problem von Shephard \"uber die projektionen
konvexer k\"orper}, Math. Z. {\bf 101} (1967), 71--82.

\bibitem[S2]{S2} R.~Schneider, {\it Convex bodies: The Brunn-Minkowski theory,}
Cambridge Univ. Press, Cambridge, 1993.

\bibitem[Sh]{Sh} G.~C.~Shephard, {\it Shadow sysytems of convex bodies}, Israel J. Math. {\bf 2} (1964)'
229--306.

\end{thebibliography}
\end{document}